\documentclass[11pt,titlepage]{article}
\usepackage{amsfonts,amssymb}
\newtheorem{theorem}{Theorem}
\newtheorem{lemma}{Lemma}

\newtheorem{definition}{Definition}
\title{Hamiltonian Systems of Negative Curvature are Hyperbolic}
\author{A.~A.~Agrachev
\and N.~N.~Chtcherbakova}
\date{39/2004/M}
\begin{document}
\maketitle
\begin{abstract}
The {\it curvature} and the {\it reduced curvature} are basic
differential invariants of the pair: $\langle$Hamiltonian system, Lagrange
distribution$\rangle$ on the symplectic manifold. We show that negativity
of the curvature implies that any bounded semi-trajectory of the
Hamiltonian system tends to a hyperbolic equilibrium, while negativity
of the reduced curvature implies the hyperbolicity of any compact
invariant set of the Hamiltonian flow restricted to a prescribed
energy level. Last statement generalizes a well-known property
of the geodesic flows of Riemannian manifolds with negative
sectional curvatures.
\end{abstract}
\subsubsection*{1 Regularity and Monotonicity} Smooth objects are supposed to be
$C^\infty$ in this note; the results remain valid for the class
$C^k$ with a finite and not large $k$ but we prefer not to specify
the minimal possible $k$.

Let $M$ be a $2n$-dimensional symplectic manifold endowed with a
symplectic form $\sigma$. A {\it Lagrange distribution}
$\Delta\subset TM$ is a smooth vector subbundle of $TM$ such that
each fiber $\Delta_x=\Delta\cap T_xM,\ x\in M,$ is a Lagrange
subspace of the symplectic space $T_xM$; in other words,
$\dim\Delta_x=n$ and $\sigma_x(\xi,\eta)=0\
\forall\xi,\eta\in\Delta_x$.

Basic examples are cotangent bundles endowed with the standard
symplectic structure and the ``vertical" distribution:
$$
M=T^*N,\ \Delta_x=T_x(T^*_qN),\quad \forall x\in T^*_qN,\ q\in N.
\eqno (1)
$$

Let $h\in C^\infty(M)$; then $\vec h\in\mbox{Vec}M$ is the
associated to $h$ Hamiltonian vector field:
$dh=\sigma(\cdot,\vec h)$. We assume that $\vec h$ is a complete
vector field, i.e. solutions of the Hamiltonian system
$\dot x=\vec h(x)$ are defined on the whole time axis. We may
assume that without a lack of generality since we are going to
study dynamics of the Hamiltonian system on compact subsets of $M$
and may reduce the general case to the complete one by the usual
cut-off procedure.

The generated by $\vec h$ Hamiltonian flow is denoted by $e^{t\vec
h},\ t\in\mathbb R$. Other notations:
$\bar\Delta\subset\mbox{Vec}M$ is the space of sections of the
Lagrange distribution $\Delta$; $[v_1,v_2]\in\mbox{Vec}M$ is the
Lie bracket (the commutator) of the fields
$v_1,v_2\in\mbox{Vec}M,\ [v_1,v_2]=v_1\circ v_2-v_2\circ v_1$.

\begin{definition} We say that $\vec h$ is regular at $x\in M$
with respect to
the Lagrange distribution $\Delta$ if
$\{[\vec h,v](x):v\in\bar\Delta\}=T_xM$.
\end{definition}

An effective version of Definition 1 is as follows:
Let $v_i\in\bar\Delta,\ i=1,\ldots,n$ be such that the vectors
$v_1(x),\ldots,v_n(x)$ form a basis of $\Delta_x$; then $\vec h$
is regular at $x$with respect to $\Delta$ if and only if the
vectors
$$
v_1(x),\ldots,v_n(x), [\vec h,v_1](x),\ldots,[\vec h,v_n](x)
$$ form a
basis of $T_xM$.

We define a bilinear mapping
$g^h:\bar\Delta\times\bar\Delta\to C^\infty(M)$ by the formula:
$$
g^h(v_1,v_2)=\sigma([\vec h,v_1],v_2).
$$
\begin{lemma} $g^h(v_2,v_1)=g^h(v_1,v_2),\ \forall
v_1,v_2\in\bar\Delta$ and $g^h(v_1,v_2)(x)$ depends only on
$v_1(x),v_2(x)$.
\end{lemma}
{\bf Proof.} Hamiltonian flows preserve $\sigma$ and $\sigma$
vanishes on $\bar\Delta$. Using these facts, we obtain:
$$
0=\sigma(v_1,v_2)=\left(e^{t\vec h*}\sigma\right)(v_1,v_2)=
\sigma(e^{t\vec h}_*v_1,e^{t\vec h}_*v_2).
$$
Differentiation of the identity
$0=\sigma(e^{t\vec h}_*v_1,e^{t\vec h}_*v_2)$
with respect to $t$ at $t=0$ gives:
 $0=\sigma([\vec h,v_1],v_2)+\sigma(v_1,[\vec h,v_2])$. Now the
anti-symmetry of $\sigma$ implies the symmetry of $g^h$. Moreover,
$g^h$ is $C^\infty(M)$-linear with respect to each argument, hence
$g^h(v_1,v_2)(x)$ depends only on $v_1(x),v_2(x). \quad \square$

Let $x\in M,\ \xi_i\in\Delta_x,\ \xi_i=v_i(x),\ v_i\in\Delta,\
i=1,2$. We set $g^h_x(\xi_1,\xi_2)=g^h(v_1,v_2)(x)$. According to
Lemma 1, $g^h_x$ is a well-defined symmetric bilinear form on
$\Delta_x$. It is easy to see that the regularity of $h$ at $x$ is
equivalent to the nondegeneracy of $g^h_x$.

If $M=T^*N$ and $\Delta$ is the vertical distribution (see (1)),
\linebreak then $g^h_x=D^2_x(h|_{T^*_qN})$, where $x\in T^*_qN$. The last
equation can be easily checked in local coordinates. Indeed, local
coordinates defined on a neighborhood $O\subset N$ provide the
identification of $T^*N\bigr|_O$ with $\mathbb R^n\times\mathbb
R^n=\{(p,q):p,q\in\mathbb R^n\}$ such that $T^*_qN$ is identified
with $\mathbb R^n\times\{q\}$, the form $\sigma$ is identified
with $\sum\limits_{i=1}^ndp_i\wedge dq_i$ and the field $\vec h$
with $\sum\limits_{i=1}^n\left(\frac{\partial h}{\partial p_i}
\frac{\partial}{\partial q_i}-\frac{\partial h}{\partial q_i}
\frac{\partial}{\partial p_i}\right)$. The fields
$\frac{\partial}{\partial p_i}$ form a basis of the vertical
distribution and
$
g^h\left(\frac{\partial}{\partial p_i},\frac{\partial}{\partial
p_j}\right)=-\left\langle dq_j,\left[\sum\limits_{i=1}^n
\left(\frac{\partial h}{\partial p_i}
\frac{\partial}{\partial q_i}-\frac{\partial h}{\partial q_i}
\frac{\partial}{\partial p_i}\right),\frac{\partial}{\partial
p_i}\right]\right\rangle=\frac{\partial^2h}{\partial p_i\partial p_j}.
$

\begin{definition} We say that a regular Hamiltonian field $\vec
h$ is monotone at $x\in M$ with respect to $\Delta$ if $g_x^h$ is a
sign-definite form.
\end{definition}

\subsubsection*{2 The Curvature}
Let $X_1,X_2$ be a pair of transversal $n$-dimensional subspaces
of $T_xM$, then $T_xM=X_1\oplus X_2$. We denote by
$\pi_x(X_1,X_2)$ the projector of $T_xM$ on $X_2$ parallel to
$X_1$. In other words, $\pi_x(X_1,X_2)$ is a linear operator
characterized by the relations $\pi_x(X_1,X_2)\bigr|_{X_1}=0,\
\pi_x(X_1,X_2)\bigr|_{X_2}=\mathbf 1$.

Now consider the family of subspaces
$J_x(t)=e^{-t\vec h}_*\Delta_{e^{t\vec h}(x)}\subset T_xM$, where
$\vec h$ is a regular Hamiltonian field; in particular,
$J_x(0)=\Delta_x$. It is easy to check that the regularity of
$\vec h$ implies the transversality of $J_x(t)$ and $J_x(\tau)$
for $t\ne\tau$, if $t$ and $\tau$ are close enough to 0. Hence
$\pi_x(J_x(t),J_x(\tau))$ is well-defined and smooth with respect
to $(t,\tau)$ in a neighborhood of $(0,0)$ with the removed
diagonal $t=\tau$. The mapping
$(t,\tau)\mapsto\pi_x(J_x(t),J_x(\tau))$ has a singularity at the
diagonal, but this singularity can be controlled. In particular,
the following statement is valid:

\begin{lemma}(see \cite{ag}). For any regular field $\vec h$,
$$
\frac{\partial^2}{\partial t\partial\tau}
\left(\pi_x(J_x(t),J_x(\tau))\bigr|_{\Delta_x}\right)\Bigr|_{\tau=0}=
t^{-2}\mathbf 1+R^h_x+O(t)\quad\mathrm{as}\ t\to 0,
$$
where $R^h_x\in\mathrm{gl}(\Delta_x)$ is a self-adjoint operator
with respect to the scalar product $g^h_x$, i.e.
$g^h_x(R^h_x\xi_1,\xi_2)=g^h_x(\xi_1,R^h_x\xi_2),\
\forall\xi_1,\xi_2\in\Delta_x$.
\end{lemma}
We set $r^h_x(\xi)=g^h_x(R^h_x\xi,\xi)$.
\begin{definition} Operator $R^h_x$ and quadratic form $r_x^h$
are called the curvature operator and the curvature form of $\vec h$
at $x$ with respect to $\Delta$. We say that $\vec h$ has a
negative (positive) curvature at $x$ if
$r^h_x(\xi)g^h_x(\xi,\xi)<0$ ($>0$), $\forall
\xi\in\Delta_x\setminus\{0\}$.
\end{definition}

It follows from the definition that only monotone fields may have
negative or positive curvature. If $\vec h$ is monotone at $x$,
then $R^h_x$ has only real eigenvalues and negativity (positivity)
of the curvature is equivalent to the negativity (positivity) of
all eigenvalues of $R^h_x$.

Let us give a coordinate presentation of $R^h_x$. Fix local
coordinates $(p,q),\ p,q\in\mathbb R^n$ in a neighborhood of $x$
in $M$ in such a way that $\Delta_x\cong\{(p,0):p\in\mathbb
R^n\}$. Let $(p(t;p_0),q(t;p_0))$ be the trajectory of the field
$\vec h$ with the initial conditions $p(0;q_0)=p_0,\ q(0;p_0)=0$.
We set $S_t=\frac{\partial q(t;p_0)}{\partial p_0}\bigr|_{p_0=0}$;
regularity of $\vec h$ is equivalent to the nondegeneracy of the
$n\times n$-matrix $\dot S_0=\frac{dS_t}{dt}\bigr|_{t=0}$. The
curvature operator is presented by the matrix Schwartzian
derivative:
$$
R^h_x=1/2\dot S_0^{-1}\stackrel{\ldots}{S}_0-3/4(\dot
S_0^{-1}\ddot S_0)^2.
$$
Examples:\begin{enumerate}\item {\sl Natural mechanical system},
$M=\mathbb R^n\times\mathbb R^n,\
\sigma=\sum\limits_{i=1}^ndp_i\wedge dq_i,\
\Delta_{(p,q)}=(\mathbb R^n,0),\
h(p,q)=1/2\|p|^2+U(q)$; then $R^h_{(p,q)}=\frac{d^2U}{dq^2}$.
\item {\sl Riemannian geodesic flow}, $M=T^*N$ and $h\bigr|_{T^*_qN}$
is a positive quadratic form $\forall q\in N$; then $h$ is
actually a Riemannian structure on $N$ which identifies the tangent and
cotangent bundles and we have: $R^h_x\xi=\mathcal R(x',\xi')x'$,
where $\mathcal R$ is the Riemanian curvature tensor and
$x',\xi'\in T_qM$ are obtained from $x,\xi\in T^*_qM$ by the
``raising of the indices".
\item {\sl Mechanical system on a Riemannian manifold}, $M=T^*N$
and $h$ is the sum of the Riemannian Hamiltonian from Example 2
and the function $U\circ\pi$, where $\pi:T^*N\to N$ is standard
projection and $U$ is a smooth function on $N$. Then
$R^h_x\xi=\mathcal R(x',\xi')x'+\nabla_{\xi'}(\nabla U)$, where
$\nabla_{\xi'}$ is the Riemannian covariant derivative.
\end{enumerate}

Now we introduce a {\it reduced curvature form} $\hat r^h_x$
defined on $\Delta_x\cap\ker d_xh$ and related to the restriction
of the Hamiltonian system on the prescribed energy level. To do
that, we need some notations. Symplectic form $\sigma_x$ on $T_xM$
induces a nondegenerate pairing of $\Delta_x$ and $T_xM/\Delta_x$.
Hence there exists a unique linear mapping $G_x:\Delta_x\to
T_xM/\Delta_x$ such that
$g_x(\xi_1,\xi_2)=\sigma_x(G_x\xi_1,\xi_2),\
\forall\xi_1,\xi_2\in\Delta_x$. The mapping $G_x$ is invertible
since the form $g_x$ is nondegenerate. Let $\Pi_x:T_xM\to
T_xM/\Delta_x$ be the canonical projection. We set
$v(x)=G^{-1}_x\Pi_x\vec h(x)$; then $v$ is a smooth section of
$\Delta$, i.e. $v\in\bar\Delta$.

Assume that $\vec h$ is a monotone field and $\vec h(x)\notin\Delta_x$;
the reduced curvature
form is defined by the formula:
$$
\hat r^h_x(\xi)=r^h_x(\xi)+\frac{3\sigma_x([\vec h,[\vec
h,v]](x),\xi)^2}{4g_x(v(x),v(x))},\quad \xi\in\Delta_x\cap\ker d_xh.
$$
In Ex.~1, we obtain:
$\hat r^h_{(p,q)}(\xi)=r^h_{(p,q)}(\xi)+\frac
3{|p|^2}\langle\frac{dU}{dq},\xi\rangle^2$.
In Ex.~2, $\hat r^h_x(\xi)=r^h_x(\xi)$. Finally, in Ex.~3
(which includes both Ex.~1 and Ex.~2) we have:
$\hat
r^h_x(\xi)=r^h_x(\xi)+\frac{3g_x(d_qU,\xi)^2}{2(h(x)-U(q))}$, where
$q=\pi(x)$.

We say that $\vec h$ has a {\it negative (positive) reduced
curvature at} $x$ if \linebreak
$\hat r^h_x(\xi)g^h_x(\xi,\xi)<0$ ($>0$), $\forall
\xi\in\Delta_x\cap\ker d_xh\setminus\{0\}$.

\subsubsection*{3 Main Results}

\begin{theorem} Let $\vec h$ be a monotone field and $x_0\in M$.
Assume that the semi-trajectory $\{e^{t\vec h}(x_0):t\ge 0\}$ has
a compact closure and $\vec h$ has a negative curvature at each
point of its closure. Then there exists
$x_\infty=\lim\limits_{t\to+\infty}e^{t\vec h}(x_0)$, where $\vec
h(x_\infty)=0$ and $D_{x_\infty}\vec h$ is hyperbolic (i.e. $D_{x_\infty}\vec h$
has no eigenvalues on the imaginary axis).
\end{theorem}

\noindent{\bf Remark.} Monotonicity of $\vec h$ is equivalent to
the monotonicity of $-\vec h$ and $R^{-h}_x=R^h_x$; hence Theorem~1
can be applied to the negative time semi-trajectories of the field
$\vec h$ as well.

\medskip\noindent Example. Consider a natural mechanical system
(Ex.~1 in Sec.~2) where $U(q)$ is a strongly concave function, then any
bounded semi-trajectory of $\vec h$ satisfies conditions of
Theorem 1.

\begin{theorem} Let ${\vec h}$ be a monotone field, $S$ be a
compact invariant subset of the flow $e^{t\vec h}$ contained in
a fix level set of $h$, $S\subset h^{-1}(c)$, and
$\vec h(x)\notin\Delta_x\
\forall x\in S$. If $\vec h$ has a negative reduced curvature at
each point of $S$, then $S$ is a hyperbolic set of the flow
$e^{t\vec h}\Bigr|_{h^{-1}(c)}$ (see \cite[Sec.~17.4]{kh} for the
definition of a hyperbolic set).
\end{theorem}

\noindent Example. Mechanical system on a Riemannian manifold
(Ex.~3 in Sec.~2). Let $\kappa_q$ be the maximal sectional
curvature of the Riemannian manifold $N$ at $q\in N$. Then any
compact invariant set $S$ of the flow $e^{t\vec h}\Bigr|_{h^{-1}(c)}$
such that the projection of $S$ to $N$ is contained in the domain
$$
\left\{q\in N:\kappa_q<0,\
\|\nabla^2_qU\|+\left(\frac 3{2(c-U(q))}+|\kappa_q|\right)\|\nabla_qU\|^2
<2|\kappa_q|(c-U(q))\right\}
$$
is hyperbolic. In particular, if $N$ is a compact Riemannian
manifold of a negative sectional curvature, then
$e^{t\vec h}\Bigr|_{h^{-1}(c)}$ is an Anosov flow for any big
enough $c$. Last statement generalizes a classical result on
geodesic flows.

\medskip Both theorems are based on the structural equations
derived in \cite{ag}. These equations are similar to the
standard linear differential equation for Jacobi vector fields in
Riemannian Geometry with the curvature operators $R^h_x$ playing
the same role as the Riemannian curvature. In particular, the
proof of Theorem~2 simply simulates the proof of the correspondent
classical result on geodesic flows. Theorem~1 describes a new
phenomenon, which is not performed by geodesic flows. Indeed, if
the curvature is negative, then the operators $R^h_x$ are
nondegenerate, while in the Riemannian case (Ex.~2 in Sec.~2) we
have $R^h_xe(x)=0$, where $e$ is the Euler field (i.e. the field
generating homothety of the fibers $T^*_qN$).


\begin{thebibliography}{9}
\bibitem{ag} A. A. Agrachev, R. V. Gamkrelidze, {\it Vector fields on
$n$-foliated $2n$-dimensional manifolds}. J. Mathematical
sciences, to appear
\bibitem{kh} A. B. Katok, B. Hasselblatt, {\it, Introduction to
the modern theory of dynamical systems}. Cambridge Univ. Press,
1997
\end{thebibliography}
\end{document}